\documentclass{elsart}

\usepackage{times}
\usepackage{amssymb}
\usepackage{amsmath}

\newcommand{\ignore}[1]{}

\newcommand{\dd}{\mathrm{d}}
\newcommand{\ii}{\mathrm{i}}
\newcommand{\R}{\mathbb{R}}

\newcommand{\Z}{\mathbb{Z}}

\newcommand{\fourth}{\mbox{$\frac{1}{4}$}}
\newcommand{\pdpd}[2]{\frac{\partial #1}{\partial #2}}
\newcommand{\dede}[2]{\frac{{\mathrm{d}} #1}{{\mathrm{d}} #2}}

\newcommand{\sgn}{\operatorname{sgn}}

\newcommand{\order}{\operatorname{order}}
\newcommand{\Op}{\operatorname{Op}}

\newcommand{\WF}{\operatorname{WF}}

\newcommand{\angleset}{I}
\newcommand{\rangset}{I'}

\newcommand{\Arho}{{A_\rho}}
\newcommand{\arho}{{a_\rho}}

\begin{document}
\begin{frontmatter}
\title{A pseudodifferential equation with damping for
one-way wave propagation in inhomogeneous acoustic media}
\author{Christiaan C. Stolk}
\address{\hskip-\parindent
        Christiaan C. Stolk,
        Centre de Math\'{e}matiques,
        Ecole Polytechnique,
        91128 Palaiseau Cedex,
        France, email: stolk@math.polytechnique.fr}
\date{December, 2003.}

\begin{abstract}
A one-way wave equation is an evolution equation in one of the space
directions that describes (approximately) a wave field. The exact wave
field is approximated in a high frequency, microlocal sense. Here we
derive the pseudodifferential one-way wave equation for an
inhomogeneous acoustic medium using a known factorization argument. We
give explicitly the two highest order terms, that are necessary for
approximating the solution. A wave front (singularity) whose
propagation velocity has non-zero component in the special direction
is correctly described. The equation can't describe singularities
propagating along turning rays, i.e.\ rays along which the velocity
component in the special direction changes sign. We show that
incorrectly propagated singularities are suppressed if a suitable
dissipative term is added to the equation.

\noindent
\raisebox{-2mm}[3mm][0mm]{{\em AMS Subject Classification:} 35L05, 35S10}
\end{abstract}



\end{frontmatter}

\section{Introduction}

In this paper we consider one-way wave equations for inhomogeneous
acoustic media in $n$ dimensions, where $n \geq 2$. We assume there is
a special space direction, which we call depth or the vertical
direction, with coordinate denoted by $z$, the other directions are
called lateral or horizontal and are denoted by $x$. The time
coordinate is denoted by $t$. The medium is described by its slowness
(inverse velocity) $\nu = \nu(z,x)$ and its mass density $\rho =
\rho(z,x)$.  Let $U=U(z,x,t)$ denote the acoustic wave field, and
$F=F(z,x,t)$ a volume source, then the acoustic equation is given by
\begin{equation}     \label{eq:acoustic_second_order}
  \left( - \rho^{-1} \nu^2 \partial_t^2
    + \sum_{j=1}^{n-1} \partial_{x_j} \rho^{-1} \partial_{x_j}
    + \partial_z \rho^{-1} \partial_z \right) U
  = F .
\end{equation}

A one-way wave equation is an equation describing only downward
propagating waves, whose propagation velocity has positive vertical
component, or only upward propagating waves, with negative vertical
component of propagation velocity. As we will discuss here, by solving
a one-way equation by progressing in depth in one direction an
approximation to (part of the) wave field is obtained.  One-way
equations are used in applications in geophysics (see e.g.\
\cite{claerbout85}) and ocean acoustics, for which many different
numerical methods have been developed. References and a discussion of
some of these can be found in e.g.\
\cite{HalpernTrefethen1988,DeHoopLeRousseauWu2000}.

In an inhomogeneous medium reflections occur for waves with wave
length comparable to the scale of the medium variations. Only the
high-frequency part can be expected to be computed using a method
progressing in one direction only. However, this is not necessarily a
disadvantage. In seismic imaging and migration an asymptotic limit is
implicit, and the absence of reflections is often an advantage.

The high-frequency (or singular) part of solutions to wave equations
is very well understood. High frequency waves propagate along rays,
curved trajectories in space. We use the theory of microlocal analysis
about this (see e.g.\ the books
\cite{duistermaat96,hormander85a,taylor81}).  Thus we consider
approximation of the wave field modulo an error that is $C^\infty$,
hence in Fourier space goes to zero when the frequency goes to
infinity, faster than any negative power of the frequency.  The
singularities of the function, that is the part that does not go to
zero rapidly when the frequency becomes large, can be localized in
position and direction using the wave front set of H\"ormander (see
the mentioned references or \cite{hormander83a}). For propagating
singularities the position and direction determine the ray. Using
these ideas we describe precisely how the wave field is approximated
by solving a one-way wave equation.

Singularities propagating with velocity that is not horizontal can be
described by a pseudodifferential evolution equation in $z$.  Such
equations are obtained from a factorization or decoupling argument,
see e.g.\ \cite{taylor75} or the similar treatment in \cite[section
9.1]{taylor81}, or \cite[section 23.2]{hormander85a}. However, such
equations are not defined at the point where the ray is tangent to
horizontal. Moreover, it is clear that waves propagated along turning
rays cannot be computed by progressing in one direction in depth
only. In practice it is desirable that such wave fronts are
suppressed.

In this paper we give a pseudodifferential evolution equation in $z$
that includes such a suppression. Thus we obtain a model for numerical
one-way wave equation methods. We give sufficient detail so that the
equation can in principle be used as a starting point for
discretization.

Our one-way wave equation is obtained in two steps. First we extend
the
factorization argument of Taylor \cite{taylor75}. We compute the
pseudodifferential equation referred to above to the two highest
orders that are at least needed for a highest order approximation of
the solution (for such a pseudodifferential equation, as for the
ordinary wave equation, the leading order term describes the behavior
of the rays, while the next term describes the amplitudes to highest
order).  We also allow for a normalization in the definition of the
down- and upward propagating parts of the wave field. With a suitably
chosen normalization the equation becomes unitary, microlocally. We
then modify this equation in order to suppress any wave fronts
propagating along turning rays. We show that the solutions to the
one-way equation approximate microlocally the real solution.

The precise formulation of these results is the subject of the next
section. Sections~\ref{sec:decoupling}
and~\ref{sec:prove_extrapolation} contain the proofs of two theorems.

\section{Pseudodifferential one-way wave equation and approximation of
the solutions}

We first introduce some notation.
The Fourier variables corresponding to $x,z$ and $t$ will
be denoted by $\xi$, $\zeta$ and $\tau$, with Fourier transform defined by
\begin{equation}
  \widehat{f}(\xi,\tau) =
  \int \int \e^{-\ii(\xi \cdot x + \tau t)} f(x,t) \, \dd x \, \dd t .
\end{equation}
A function
$\psi = \psi(x,t,\xi,\tau)$ in $C^{\infty}(\R^n \times \R^n)$ is a symbol of
order $m$, if there are constants $C_{\alpha,\beta}$ such that
\begin{equation} \label{eq:standard_symbol_estimate}
  | \partial_{x,t}^\alpha \partial_{\xi,\tau}^\beta \psi(x,t,\xi,\tau) | <
    C_{\alpha,\beta}(1+\|(\xi,\tau)\|)^{m-|\beta|}
\end{equation}
(see a text on pseudodifferential operators).
Here $\alpha,\beta$ are a multi-indices $\alpha = (\alpha_1, \ldots, \alpha_n)$,
$\alpha_j \in \Z_{\geq 0}$, $|\alpha| = \alpha_1 + \ldots + \alpha_n$.
Associated with such a symbol $\psi$ is a pseudodifferential operator, that
will be denoted by $\psi(x,t,D_x,D_t)$, and is given by
\begin{equation}
  \psi(x,t,D_x,D_t) f = (2\pi)^{-n} \int
    \psi(x,t,\xi,\tau) \widehat{f}(\xi,\tau)
    \e^{\ii (\xi \cdot x + \tau t)} \, \dd \xi \, \dd \tau .
\end{equation}
The set of symbols is denoted by $S^m(\R^n \times \R^n)$.
We will typically encounter operators acting in $(x,t)$, depending on
$z$, but independent of $t$ (of convolution type in $t$), i.e.\ with
symbols $\psi(z,x,\xi,\tau)$.
A symbol is in $S^{-\infty}$ if for any $N$ there is a constant $C$
such that $| \partial_{x,t}^\alpha \partial_{\xi,\tau}^\beta \psi |
< C(1+\|(\xi,\tau)\|)^{-N}$
The corresponding operator then maps distributions to $C^\infty$.  We
will write $\psi_1 \sim \psi_2$ if $\psi_1 - \psi_2 \in
S^{-\infty}$.

We recall that, if $f(y)$, $y \in \R^n$ is a distribution, H\"ormander's
wave front set $\WF(f)$ contains points $(y,\eta)$, contained in the cotangent
space or phase space $\R^n \times \R^n$,
that give positions and directions
associated with the singularities of $f$, see \cite[section
1.3]{duistermaat96}, or \cite[section 8.1]{hormander83a}. This
set is conic, i.e.\ if $(y,\eta) \in \WF(f)$ then $(y,\lambda \eta)
\in \WF(f)$ for all $\lambda >0$. If $\Gamma$ is a conic subset of
$\R^n \times \R^n$ we say that $f \equiv g$ microlocally on $\Gamma$
if $\WF(f - g) \cap \Gamma = \emptyset$.

Denote by $p(z,x,\zeta,\xi,\tau) = \rho(z,x)^{-1} \nu(z,x)^2 \tau^2
  - \rho^{-1} (\xi^2 + \zeta^2)$ the principal symbol of $P$.
If $U$ satisfies $P U = 0$, then the singularities of $U$ are in
the characteristic set given by
\begin{equation} \label{eq:characteristic}
  p(z,x,\zeta,\xi,\tau) = 0 .
\end{equation}
They propagate along null bicharacteristics, curves in the cotangent
space contained in the set given by (\ref{eq:characteristic})
that are solutions to the Hamilton vector field of $p$ (see
\cite[Theorem 6.2.1]{taylor81} or \cite[section 23.1]{hormander85a}).
Parameterizing by time, the differential equations for the
null bicharacteristics are
\begin{align}
  \dede{x}{t} = {}& - \nu(z,x)^{-2} \tau^{-1} \xi , &
    \dede{z}{t} = {}& - \nu(z,x)^{-2} \tau^{-1} \zeta , \\
  \dede{\xi}{t} = {}& - \tau \nu^{-1}\pdpd{\nu}{x} , &
    \dede{\zeta}{t} = {}& - \tau \nu^{-1}\pdpd{\nu}{z} .
\end{align}
So singularities with $-\frac{\zeta}{\tau} > 0$, satisfy $\dede{z}{t} >
0$ (downgoing) and waves with $-\frac{\zeta}{\tau} < 0$ satisfy
$\dede{z}{t} < 0$ (upgoing).

If a point $(z,x,\xi,\tau)$ is given with $\|\xi\| < \nu(z,x) |\tau|$,
then there are two solutions $\zeta$ to (\ref{eq:characteristic}).
These will be denote by $\pm b = \pm b(z,x,\xi,\tau)$, where
\begin{equation} \label{eq:define_b_pm}
  b = -\tau \nu \sqrt{ 1 - \tau^{-2} \nu^{-2} \xi^2 } .
\end{equation}
The sign is such that $\pm b$ corresponds to rays with $\pm
\pdpd{z}{t}>0$.
We define a set of points $(z,x,\xi,\tau)$ associated with propagation
angles with the vertical $< \theta$, by
\begin{equation} \label{eq:set_regular}
  \angleset'_{\theta} = \{ (z,x,\xi,\tau) \, | \, \tau\neq 0
    \text{ and }
    \| \nu(z,x)^{-1} \tau^{-1} \xi \| \leq \sin(\theta) \} .
\end{equation}

To obtain the microlocal one-way wave equations we follow the
factorization argument of Taylor \cite{taylor75}. In this
factorization it is assumed that singularities of $F$ and therefore
also those of $U$ are not in the set given by $\xi=\tau=0, \zeta \neq
0$. By (\ref{eq:characteristic}) there is a constant $C$ (assuming
that $\nu$ is bounded) such that for the propagating singularities we
have
\begin{equation} \label{eq:max_k_z}
  |\zeta| < C |\tau| .
\end{equation}
The decoupling will be done microlocally where the propagation angle
$\theta$ is smaller than some fixed angle $\theta_2$. We define a
subset of phase space $\R^{n+1} \times \R^{n+1}$ associated to such
angles by
\begin{equation}
  \angleset_{\theta_2} =  \{ (z,x,t,\zeta,\xi,\tau) \,|\,
  (z,x,\xi,\tau) \in \rangset_{\theta_2}, |\zeta| < C |\tau| \} .
\end{equation}
Let $I_{\pm,\theta_2}$ denote the subsets of $\angleset_{\theta_2}$ with
$\mp \tau^{-1} \zeta > 0$.

The down- and upgoing components $(u_+,u_-)$ are obtained from
$(U,\rho^{-1} \pdpd{U}{z})$ by a $2\times2$ pseudodifferential matrix
$Q = Q(z,x,D_x,D_t)$. We assume that $Q$ is elliptic, i.e.\ there is a
matrix pseudodifferential operator $W$ such that $W Q \sim Q W
\sim \operatorname{Id}$.  We take the liberty to denote $Q^{-1} = W$,
even though this is not quite correct.
We also define sources $(f_+,f_-)$. We have
\begin{align} \label{eq:define_u_pm_first}
  \left(\begin{matrix} u_+\\u_- \end{matrix} \right)
  = {}& Q^{-1} \left( \begin{matrix} U \\ \rho^{-1} \pdpd{U}{z}
        \end{matrix} \right) , &
  \left(\begin{matrix} f_+\\f_- \end{matrix} \right)
  = {}& Q^{-1} \left( \begin{matrix} 0 \\ F
        \end{matrix} \right) .
\end{align}
We have the following result about one-way wave equations for
$(u_+,u_-)$. The proof is the subject of section~\ref{sec:decoupling}.

\begin{thm} \label{th:decoupling}
For suitably chosen $Q$ and $B_\pm = B_\pm(z,x,D_x,D_t)$ the equation
\begin{equation} \label{eq:PU_F_microlocally}
  P U \equiv F \text{ microlocally on $\angleset_{\theta_2}$} .
\end{equation}
holds if and only if
\begin{align} \label{eq:one_way_microlocal}
  P_{0,+} u_+ \stackrel{\rm def}{=}
    ( \partial_z - \ii B_+(z,x,D_x,D_t) ) u_+ \equiv {}& f_+ ,
    \text{ microlocally on $\angleset_{\theta_2}$} , \text{ and } \\
  P_{0,-} u_- \stackrel{\rm def}{=}
    ( \partial_z - \ii B_-(z,x,D_x,D_t) ) u_- \equiv {}& f_- ,
    \text{ microlocally on $\angleset_{\theta_2}$} .
\end{align}
The operator $B_\pm$ can be chosen selfadjoint, with $B_\pm$ and $Q$ satisfying
\begin{align} \label{eq:symbol_B}
  B_{\pm}(z,x,\xi,\tau) = {}&
    \pm \left( b + \half \ii b^{-1} \sum_{j=1}^{n-1}
        \pdpd{b}{\xi_j} \pdpd{b}{x_j} \right)
    + \order(-1) ,
    \text{ on ${\rangset_{\theta_2}}$} , \\
  \label{eq:symbol_Q}
  Q(z,x,\xi,\tau) = {}&
     \left( \begin{matrix}
        \rho^{\frac{1}{2}} a^{-\frac{1}{4}}
        & \rho^{\frac{1}{2}} a^{-\frac{1}{4}} \\
        \ii \sgn(\tau) \rho^{-\frac{1}{2}} a^{\frac{1}{4}}
        & -\ii \sgn(\tau) \rho^{-\frac{1}{2}} a^{\frac{1}{4}}
    \end{matrix} \right)
    + \order \left(\begin{matrix} -\frac{3}{2} & -\frac{3}{2} \\
        -\frac{1}{2} & -\frac{1}{2} \end{matrix} \right) ,
  \text{ on } \rangset_{\theta_2} ,
\end{align}
where $a$ is defined by $a(z,x,\xi,\tau) =
\nu(z,x)^2 \tau^2 - \| \xi \|^2$.
If we choose $Q$ such that $U = u_+ + u_-$, that is $Q_{1,1} = Q_{1,2} =
1$, then $B_\pm$ satisfies
\begin{multline} \label{eq:B_pm_non_sa}
  B_\pm(z,x,\xi,\tau) =
    \pm \left( b + \half \ii b^{-1} \sum_{j=1}^{n-1}
        \pdpd{b}{\xi_j} \pdpd{b}{x_j} \right)
    + \fourth \ii \pdpd{a}{z} a^{-1}
    - \half \ii \pdpd{\rho}{z} \rho^{-1} \\
  + \order(-1) ,
  \text{ on $\rangset_{\theta_2}$} .
\end{multline}
\end{thm}

The highest order term in equation (\ref{eq:one_way_microlocal})
determines the Hamilton flow of the singularities, hence the rays. The
zeroeth order term determines the amplitude. Thus both terms need to
be incorporated for an accurate highest order approximation of the
solutions using (\ref{eq:one_way_microlocal}).
There are two kinds of zeroeth order terms for $B_\pm$. First the term
$\half \ii b^{-1} \sum_{j=1}^{n-1} \pdpd{b}{\xi_j}
\pdpd{b}{x_j}$. When computed explicitly it is equal to
$\half \ii \nu \xi \cdot \pdpd{\nu}{x} \tau^{-1}
(\nu^2-\tau^{-2} \xi^2)^{-3/2}$.
This term makes the operator $b + \half \ii
b^{-1} \sum_{j=1}^{n-1} \pdpd{b}{\xi_j} \pdpd{b}{x_j}$
selfadjoint up to zeroeth order. Second there is the term
$\fourth \ii\pdpd{a}{z} a^{-1} - \half \ii \pdpd{\rho}{z} \rho^{-1}$ in
(\ref{eq:B_pm_non_sa}), due to the different normalization of $Q$.

Equation (\ref{eq:one_way_microlocal}) is only microlocal.  Outside
$\rangset_{\theta_2}$ the symbol $B_\pm$ is not prescribed, but we
choose it with real principal symbol and smooth (i.e.\ without the
singularity of the square root). A wave front propagating on some
turning ray will be propagated incorrectly with this equation. To
suppress such singularities we introduce a damping term given by a
pseudodifferential operator $C = C(z,x,D_x,D_t)$. The complete one way
wave equation will be of the form
\begin{equation} \label{eq:reg_one_way}
  P_{\pm} u_\pm
    \stackrel{\rm def}{=} ( \partial_z - \ii B_{\pm}(z,x,D_x,D_t)
    + C(z,x,D_x,D_t) ) u_{\pm} = 0 .
\end{equation}
It is assumed here that a solution is sought for $z > z_0$ (for $z <
z_0$ the sign in front of $C(z,x,D_x,D_t)$ must be changed).
We assume the dissipative term
is $0$ for waves propagating with angle smaller than some
given angle $\theta_1$, $\theta_1 < \theta_2$.
The operator $C$ is also a pseudodifferential operator with
homogeneous, non-negative principal symbol $c(z,x,\xi,\tau)$. We let
the order be 1, so that the length scale associated with the
decay is proportional to wave
length, but this is not essential. Its main property will be
\begin{align} \label{eq:C_main}
  c(z,x,\xi,\tau) & = 0
     && \text{ for $(z,x,\xi,\tau) \in \rangset_{\theta_1}$}\\
  c(z,x,\xi,\tau) & \geq \eta (\xi^2 + \tau^2)^{\frac{1}{2}}
     && \text{ for $(z,x,\xi,\tau)$ outside $\rangset_{\theta_2}$} ,
\end{align}
where $\eta$ is some positive constant.
In addition there is the condition that when $C = 0$, then also a
number of its derivatives are zero, see the precise formulation below.

We consider the approximation of solutions to the equation
(\ref{eq:acoustic_second_order}) with right hand side $0$, given that
the solution $U$ has only singularities propagating in the $+$
direction (or only in the $-$ direction). So suppose that $U$ satisfies
\begin{equation} \label{eq:PU_equiv_0}
  P U \equiv 0 \text{ for $z>z_0$},
\end{equation}
and assume it has only  singularities propagating in the $+$ direction
at $z=z_0$, in other words
\begin{equation} \label{eq:WF_U_directive}
  WF(U) \cap \{ z = z_0 , \tau^{-1} \zeta > 0 \} = \emptyset .
\end{equation}
We also assume that the singularities of $PU$ at $z=z_0$ satisfy
(\ref{eq:max_k_z}) at $z=z_0$, which implies that the restriction
$U|_{z_0}$ is well defined.
Let $Q_+ = Q_{1,1}, Q_- = Q_{1,2}$. The approximate solution is then
given by $Q_+ u_+$, where $u_+$ is the solution of
\begin{align}
  u_+ \big|_{z_0} = {}& Q_+(z_0)^{-1} U|_{z_0} ,
        \label{eq:u_pm_init}\\
  P_+ u_+ = {}& 0 , & & \text{ for $z>z_0$} , \label{eq:u_pm_psdiff_eq}
\end{align}

With a point $(z_0,x,t,\xi,\tau) \in \rangset_{\theta}$ and a time $t$
there are two associated null bicharacteristics corresponding to the
two possible values of $\zeta = \pm b(z_0,x,\xi,\tau)$. They can be
parameterized by $z$-coordinate of the ray as long as the angle of the
velocity vector is smaller or equal than $\theta$. We let
$[Z_{\rm min}(z_0,x,\xi,t,\tau,\theta), \linebreak[2]
Z_{\rm max}(z_0,x,\xi,t,\tau,\theta)]$ be the maximal interval where this is
the case, and we denote the bicharacteristic by
$\gamma_\pm(z,z_0,x,t,\xi,\tau)$. Let $c(z,x,\xi,\tau)$ be the
principal symbol of $C$. It was shown in \cite{StolkMixParametrix1}
that in a solution operator for (\ref{eq:u_pm_psdiff_eq}) a
pseudodifferential factor occurs with symbol
\begin{equation} \label{eq:C_exp}
  \exp \bigg(
    - \int_{z_0}^z c(\gamma_\pm(z',z,x,t,\xi,\tau)) \, \dd z' \bigg) .
\end{equation}
This exponential is equal to $1$ when the bicharacteristic
$\gamma_\pm$ stays in $\angleset_{\theta_1}$, while it is
exponentially decaying to $0$ for $(\xi,\tau)$ to infinity if a finite
segment between $z_0$ and $z$ is outside the region where $c = 0$.
Therefore we define a subset of $\R^{n+1} \times \R^{n+1}$ that can be
reached from depth $z_0$, while staying in $\angleset_{\theta}$, by
\begin{equation}
  J_{\pm}(z_0,\theta)
  = \{ (z,x,t,\zeta,\xi,\tau) \, | \, \mp \tau^{-1} \zeta > 0 \text{ and }
        Z_{\rm min}(z,x,\xi,\tau,\theta) \leq z_0 \} .
\end{equation}
We will show $U$ can be approximated by $Q_+ u_+$ in the following way
\begin{align}
  \WF(Q_+ u_+) \subset {}& \WF(U) ,
    \label{eq:approx_prop_1}\\
  Q_+ u_+ \equiv {}& U \text{ on $J_+(z_0,\theta_1)$} ,
    \label{eq:approx_prop_2}\\
  Q_+ u_+ \equiv {}& 0 \text{ outside $J_+(z_0,\theta_2)$} .
    \label{eq:approx_prop_3}
\end{align}
The same is true with $u_+, Q_+$ and $J_+$ replaced by
$u_-,Q_-$ and $J_-$.

In \cite{StolkMixParametrix1} additional assumptions on $C$ were
made. We first give an example.  Define a scalar function $h(y)$ that
smoothly goes from constant equal to zero at $y <0$, to being positive
at $y>0$ by the formula
\begin{equation} \label{eq:h_choice}
  h(y) = \left\{ \begin{array}{ll}
        0          & y \leq 0 , \\
        \exp(-1/y)/(\exp(-1/y)+\exp(-1/(1-y))) & 0 < y < 1 , \\
        1          & y \geq 1 \end{array} \right .
\end{equation}
Now define for example
\begin{equation} \label{eq:C_example}
  C = w(z,x,\xi,\tau)
    h(\nu^{-1} \|\tau^{-1} \xi\| - \sin(\theta_1))
\end{equation}
where $w$ is homogeneous of order $1$ in $(\xi,\tau)$ and bounded
below by some constant times $\sqrt{\tau^2+\xi^2}$.

In general we make the following assumptions. We assume that $C$
is given by a sum $C = c + C^{(1)}$ where $C^{(1)}$ is of order $0$. We
will also write $C^{(0)}(z,x,\xi,\tau) = c(z,x,\xi,\tau)$. We assume that
there is an integer $L > 2$ such that the derivatives of order up to
$L$ of $C^{(0)}$, and of order $L-2$ of $C^{(1)}$ satisfy the following
bounds
\begin{multline} \label{eq:assumption_C}
  \big| \partial_z^j \partial_{x}^\alpha \partial_{\xi,\tau}^\beta
      C^{(k)}(z,x,\xi,\tau) \big|
  \leq C (1+\|(\xi,\tau)\|)^{-|\beta|-k+\frac{j+|\alpha|+|\beta|}{L}}
    (1+c(z,x,\xi,\tau))^{1-\frac{j+|\alpha|+|\beta|}{L}} , \\
    j + 2k + |\alpha| + |\beta| < L ,
\end{multline}
for some constant $C$.
It was shown in \cite{StolkMixParametrix1} that (\ref{eq:C_example})
satisfies this property for any $L$.

The initial value problems for the operators $P_{0,\pm}$ and $P_\pm$
that were defined in (\ref{eq:one_way_microlocal}) and (\ref{eq:reg_one_way})
have well defined solution operators, that we will denote by
$E_{0,\pm}(z,z_0)$ and
$E_\pm(z,z_0)$. In \cite{StolkMixParametrix1}
it was shown that these are related by a pseudodifferential operator
$K_\pm = K_\pm(z,z_0,x,D_x,D_t)$ with principal symbol (\ref{eq:C_exp}),
such that $E_\pm \sim K_\pm E_{0,\pm}$ (with $K$ in a class of symbols more
general than that given by (\ref{eq:standard_symbol_estimate})).
With the assumptions on $C$ it follows that $K$ has the property
\begin{align}
  K_\pm(z,z_0,x,\xi,\tau) \sim {}& 1 \text{ on $J_\pm(z_0,\theta_1)$} ,
    \label{eq:K_suppress_property_1}\\
  K_\pm(z,z_0,x,\xi,\tau) \in {}& S^\infty \text{ outside $J_\pm(z_0,\theta_2)$,
        if $z-z_0 > \delta$} ,
    \label{eq:K_suppress_property_2}
\end{align}
if $\delta > 0$ is some small constant. This results in the following
theorem, that gives sufficient conditions for the approximation
property of equations (\ref{eq:approx_prop_1}),
(\ref{eq:approx_prop_2}) and (\ref{eq:approx_prop_3}) to hold. The
proof is given in section~\ref{sec:prove_extrapolation}.

\begin{thm} \label{th:square_root_correct_sols}
Let $U = U(z,x,t)$ satisfy (\ref{eq:PU_equiv_0}) and
(\ref{eq:WF_U_directive}). Let $u_+$ be the solution to
(\ref{eq:u_pm_init}) and (\ref{eq:u_pm_psdiff_eq}),
where $Q_+,B_+$ are as in Theorem~\ref{th:decoupling}, and $C$ satisfies
(\ref{eq:C_main}) 
and (\ref{eq:assumption_C}).
Then there is $K$, depending on $Q_+, B_+$ and $C$, satisfying
(\ref{eq:K_suppress_property_1}) and (\ref{eq:K_suppress_property_2})
such that for $z > z_0$
\begin{equation}
  Q_+ u_+ = K U + r ,
\end{equation}
with $r \in C^\infty(]z_0,\infty[ \times \R^n)$.
\end{thm}


\section{Proof of Theorem~\ref{th:decoupling}\label{sec:decoupling}}

The computation of $Q$ and $B_\pm$ is done by writing
(\ref{eq:acoustic_second_order}) as a system of first order in $z$,
and then transforming this system.  Let $V$ be defined by $V =
\rho^{-1}\pdpd{U}{z}$, and $A, \Arho$ by
\begin{align} \label{eq:define_Arho}
  \Arho = {}& -\rho^{-1} \nu(z,x)^2 \partial_t^2
    + \sum_j \partial_{x_j} \rho^{-1} \partial_{x_j} , \\
  \label{eq:define_A}
  A = {}& A_1 = - \nu(z,x)^2 \partial_t^2 + \sum_j \partial_{x_j}^2 .
\end{align}
The principal symbols of these operators are $a(z,x,\xi,\tau)
= \linebreak[2] \nu(z,x)^2 \tau^2 - \| \xi \|^2$, \linebreak[2]
$a_\rho(z,x,\xi,\tau)
= \linebreak[2] \rho^{-1} \nu(z,x)^2 \tau^2- \rho^{-1} \| \xi \|^2$.
With these definitions, equation (\ref{eq:acoustic_second_order})
is equivalent to the following system for the vector $(U,V)$
\begin{equation} \label{eq:first_order_system}
  \pdpd{}{z} \left( \begin{matrix} U\\ V \end{matrix} \right)
  - \left( \begin{matrix} 0 & \rho \\ -A_\rho & 0 \end{matrix} \right)
    \left( \begin{matrix} U\\ V \end{matrix} \right)
  = \left( \begin{matrix} 0\\ F \end{matrix} \right) .
\end{equation}
The transformed wave field $(u_+, u_-)$ and source $(f_+,f_-)$ were
defined in (\ref{eq:define_u_pm_first}) from $(U,V)$ and $(0,F)$.
Recall that the matrix pseudodifferential operator $Q$ is elliptic,
and that $Q^{-1}$ denotes a microlocal inverse, not an exact inverse,
satisfying $Q^{-1} Q \sim Q Q^{-1} \sim \operatorname{Id}$.

A technical complication is that the operators $Q$
and $Q^{-1}$ are not pseudodifferential operators in
$(z,x,t)$ (only in $(x,t)$).
We let $\psi=\psi(z,x,D_z,D_x,D_t)$ be a microlocal cutoff
around $\xi=\tau=0$. We let its symbol $\psi(z,x,\zeta,\xi,\tau)$ be $1$
for $|\zeta| < 2C|\tau|, |\zeta| >1$ and $0$ for $|\zeta| > 3C|\tau|$.
Equation (\ref{eq:PU_F_microlocally}) is true if and only if
\begin{equation} \label{eq:reg_first_order_system}
  \psi Q^{-1} \left[
  \pdpd{}{z}
  - \left( \begin{matrix} 0 & \rho \\ -A_\rho & 0 \end{matrix} \right) \right]
    \left( \begin{matrix} U\\ V \end{matrix} \right)
  \equiv \psi Q^{-1}
    \left( \begin{matrix} 0\\ F \end{matrix} \right) ,
  \text{ microlocally on $\angleset_{\theta_2}$} .
\end{equation}
By \cite[theorem 18.1.35]{hormander85a} the operator $\psi Q^{-1}$ is
a pseudodifferential operator with symbol that equals $Q^{-1}$ modulo
$S^{-\infty}$ on $\angleset_{\theta_2}$. Using this theorem again we
can see that in (\ref{eq:reg_first_order_system})
a factor $Q Q^{-1}$ can be inserted in the left side before $(U,V)$.
So (\ref{eq:reg_first_order_system}) holds if and only if
\begin{equation} \label{eq:transformed_equation}
  \psi Q^{-1} \bigg[ \pdpd{}{z}
        - \left( \begin{matrix} 0 & \rho \\ -\Arho & 0 \end{matrix} \right)
    \bigg] Q
    \left(\begin{matrix} u_+ \\ u_- \end{matrix} \right)
  \equiv  \psi \left(\begin{matrix} f_+ \\ f_- \end{matrix} \right) ,
  \text{ microlocally on $\angleset_{\theta_2}$} .
\end{equation}
Therefore Theorem~\ref{th:decoupling} follows from the lemma that
we now state concerning the diagonalization of the operator
\begin{equation} \label{eq:transformed_operator1}
  Q^{-1} \bigg[ \pdpd{}{z}
        - \left( \begin{matrix} 0 & \rho \\ -\Arho & 0 \end{matrix} \right)
    \bigg] Q .
\end{equation}

\begin{lem} \label{lem:decoupling}
For suitably chosen $Q$ and $B$, the operator (\ref{eq:transformed_operator1})
is equal to
\begin{equation} \label{eq:almost_diagonal_operator}
  \pdpd{}{z}
    - \ii \left( \begin{matrix}  B_+ & 0 \\
                             0 & B_- \end{matrix} \right)
    + R ,
\end{equation}
where $R = R(z,x,D_x,D_t)$ is a $2 \times 2$ matrix pseudodifferential
operator of order $1$ with
symbol that is in $S^{-\infty}$ on $\rangset_{\theta_2}$.  Here we
can choose $B_\pm$ selfadjoint, with $Q,B_\pm$ satisfying (\ref{eq:symbol_Q})
and (\ref{eq:symbol_B}). We can also choose $Q$ with $Q_{1,1} = Q_{1,2} = 1$,
with $B_\pm$ satisfying (\ref{eq:B_pm_non_sa}) and $Q$ satisfying
\begin{equation}
  Q = \left( \begin{matrix} 1 & 1 \\
        -\ii \sgn(\tau) \rho^{-\frac{1}{2}} \arho^{\frac{1}{2}}
          +\order(0)
        & \ii \sgn(\tau) \rho^{-\frac{1}{2}} \arho^{\frac{1}{2}}
          +\order(0)
    \end{matrix} \right) .
\end{equation}
\end{lem}

The proof mostly follows an argument of Taylor \cite{taylor75}. Some extra
work is required to obtain the explicit expressions and the symmetry
(self-adjointness) property. It is an order by order construction resulting in
an asymptotic sum that is well defined according to standard results
(see e.g.\ \cite[Proposition 18.1.3]{hormander85a}).  The main tool is
the composition formula, which says that the product of two
pseudodifferential operators $A(y,D_y)$ and $B(y,D_y)$ is again a
pseudodifferential operator with symbol $A \# B$ given by the
asymptotic sum
\begin{equation} \label{eq:psido_composition}
  \sum_{\alpha} \frac{1}{\alpha! \, \ii^{|\alpha|}}
     \partial^\alpha_{\eta} A(y,\eta) \partial_y^\alpha B(y,\eta) .
\end{equation}
Here $\alpha!=\alpha_1! \ldots \alpha_n!$ (if $y \in \R^n$).

For the computations in the proof of Lemma~\ref{lem:decoupling}, we
use pseudodifferential operators that are microlocally the square root or
certain other powers of the operators $A$, $\Arho$. In addition we use
an operator $S(D_t)$ with symbol $-\sgn(\tau)$.  In the following
lemma we collect the needed information about these.

\begin{lem} \label{lem:square_root}
There is a pseudodifferential square root operator $\widetilde{B} =
\widetilde{B}(z,x,D_x,D_t)$, satisfying $\widetilde{B}^2 \sim A$ microlocally
on $\rangset_{\theta_2}$. Its symbol satisfies
\begin{equation} \label{eq:symbol_square_root}
  \widetilde{B}(z,x,\xi,\tau)
  = b + \half \ii b^{-1} \sum_{j=1}^{n-1}
        \pdpd{b}{\xi_j} \pdpd{b}{x_j}  + \order (-1) ,
  \text{ on $\rangset_{\theta_2}$} .
\end{equation}
There are fourth and second roots $\Arho^s$, $s=\frac{1}{4},\frac{1}{2}$,
microlocally on $\rangset_{\theta_2}$, with principal symbol that
equals $\arho(z,x,\xi,\tau)^s$ on $\rangset_{\theta_2}$. There
are inverse $\Arho^{-s}$, $s=\frac{1}{4},\frac{1}{2},1$ microlocally
on $\rangset_{\theta_2}$, with principal symbol
$\arho(z,x,\xi,\tau)^{-s}$ on $\rangset_{\theta_2}$.
The operators $\widetilde{B},\Arho^s$ can be chosen selfadjoint.
Let $S(z,x,D_x,D_t)$ be a pseudodifferential operator with symbol
equal to $-\sgn(\tau)$ on $\rangset_{\theta_2}$, $|\tau| > 1$, and selfadjoint.
Then $\widetilde{B} \sim S A^{\frac{1}{2}}$ on $\rangset_{\rm \theta_2}$.
We have
\begin{align} \label{eq:square_root_rho_equalities}
  \widetilde{B}
  = {}& S \rho^{\frac{1}{4}} \Arho^{\frac{1}{2}} \rho^{\frac{1}{4}}
        + \order(-1)   \nonumber \\
  = {}& S \rho^{-\frac{1}{4}} \Arho^{\frac{1}{4}} \rho
        \Arho^{\frac{1}{4}} \rho^{-\frac{1}{4}} + \order(-1) ,
          \text{ on $\rangset_{\theta_2}$} .
\end{align}
\end{lem}

\begin{pf}
Following the standard argument (compare e.g.\ Lemma II.6.2 in
\cite{taylor81}) we look for the square root as an asymptotic sum
$\sum_{j=0}^\infty T^{(j)}$ with each $T^{(j)}$ a pseudodifferential
operator of order $1-j$. We let $T^{(0)}$ have principal symbol $b$,
then $(T^{(0)})^2 = A + R^{(0)}$, with $R^{(0)}$ of order $1$. Now
suppose we have $T^{(j)}$, $j=0,\ldots,k$, such that
\begin{equation} \label{eq:square_root_induction}
  \left( \sum_{j=0}^k  T^{(j)}\right)^2 = A + R^{(k)}, \text{ with $R^{(k)}$ of
    order $1-k$} .
\end{equation}
Then we let $T^{(k+1)}$ have principal symbol $-\half b^{-1} R^{(k)}$,
and (\ref{eq:square_root_induction}) is valid with $k+1$ instead of
$k$ (note that this choice of $T^{(k+1)}$ is unique to highest order on
$\rangset_{\theta_2}$). In this case $R^{(0)}$ is $\ii \sum_{j=1}^{n-1}
\pdpd{b}{\xi_j} \pdpd{b}{x_j} + \order (-1)$, which leads
to (\ref{eq:symbol_square_root}). Because $b$ is real all the $T^{(j)}$
can be chosen selfadjoint.   The square root
$A_\rho^{\frac{1}{2}}$ follows similarly and its fourth root equals
the square root of the square root. The existence of microlocal
inverses is standard (see Theorem 18.1.9 in \cite{hormander85a}).

The symbol of $S$ is locally constant on $\rangset_{\theta_2}$ so $S$ commutes
microlocally with $\rho^s$ and the $\Arho^s$.
Also the commutator $[\Arho^s,\rho^{s'}]$ is a lower order operator
(of order $\frac{s}{2}-1$). To compute the square of $S \rho^{\frac{1}{4}}
\Arho^{\frac{1}{2}} \rho^{\frac{1}{4}}$ we commute a factor
$\rho^{\frac{1}{4}}$ to the left of the first $\Arho^{\frac{1}{2}}$ and
one to the right of the second $\Arho^{\frac{1}{2}}$. Using that
the multiple commutator is again an order lower, we find
\begin{equation} \label{eq:intermed}
  (S \rho^{\frac{1}{4}} \Arho^{\frac{1}{2}} \rho^{\frac{1}{4}})^2
  =  \rho^{\frac{1}{2}} \Arho \rho^{\frac{1}{2}} + \order(0)
  =  A + \order(0) , \text{ on $\rangset_{\theta_2}$.}
\end{equation}
Equation (\ref{eq:square_root_rho_equalities}) is clearly valid for
the principal symbols.
It follows from (\ref{eq:intermed}) and the fact that $T^{(1)}$ is unique
to highest order that the first equality in
(\ref{eq:square_root_rho_equalities}) is satisfied. The
second equality of (\ref{eq:square_root_rho_equalities}) follows similarly.
\qed
\end{pf}

\begin{pf*}{Proof of Lemma~\ref{lem:decoupling}}
Commuting $Q^{-1}$ and $\pdpd{}{z}$ we find that
(\ref{eq:transformed_operator1}) is equal to
\begin{equation} \label{eq:transformed_operator2}
  \pdpd{}{z}
    - Q^{-1} \left( \begin{matrix} 0 & \rho \\ -\Arho & 0 \end{matrix} \right) Q
    - \pdpd{Q^{-1}}{z} Q .
\end{equation}
We first consider the second term of (\ref{eq:transformed_operator2}), which
is a contribution of order $1$ (the third term is of order $0$). The
eigenvalues and eigenvectors of the principal symbol matrix $\left(
\begin{matrix} 0 & \rho \\ -\arho & 0 \end{matrix} \right)$
are given by
\begin{equation}
  \text{eigenvalues : } \pm \ii \rho(z,x)^{1/2} \arho(z,x,\xi,\tau)^{1/2} ,
  \text{ eigenvectors : } \left( \begin{matrix} \rho(z,x)^{1/2} \\
        \pm \ii \arho(z,x,\xi,\tau)^{1/2} \end{matrix} \right) .
\end{equation}
We first make a highest order choice $Q = Q^{(0)}$ such that the
matrix
$Q^{-1} \left( \begin{matrix} 0 & \rho \\ -\Arho & 0 \end{matrix} \right) Q$
becomes diagional on $\rangset_{\theta_2}$. Below we will add lower order
terms. We set
\begin{equation}
  Q^{(0)} = \left(\begin{matrix}
        \Arho^{-\frac{1}{4}} \rho^{\frac{1}{4}}
            & \Arho^{-\frac{1}{4}} \rho^{\frac{1}{4}} \\
        \ii S \Arho^{\frac{1}{4}} \rho^{-\frac{1}{4}}
            & - \ii S \Arho^{\frac{1}{4}} \rho^{-\frac{1}{4}}
    \end{matrix} \right) ,
    \text{ on $\rangset_{\theta_2}$} ,
\end{equation}
with $S$ as in Lemma~\ref{lem:square_root}.
Outside (\ref{eq:set_regular}) we still require that the symbol $Q^{(0)}$ is
an invertible matrix symbol of order $\left(\begin{matrix}
-\frac{1}{2}&-\frac{1}{2}\\ \frac{1}{2}&\frac{1}{2}\end{matrix}\right)$.
It is easily seen that this is possible. The inverse of $Q^{(0)}$ satisfies
\begin{equation} \label{eq:symbol_Q_inverse}
  Q^{(0)}(z,x,\xi,\tau)^{-1} = \frac{1}{2} \left(\begin{matrix}
        \rho^{-\frac{1}{4}} \Arho^{\frac{1}{4}}
            & -\ii S \rho^{\frac{1}{4}} \Arho^{-\frac{1}{4}} \\
        \rho^{-\frac{1}{4}} \Arho^{\frac{1}{4}}
            & \ii S \rho^{\frac{1}{4}} \Arho^{-\frac{1}{4}} ,
    \end{matrix}\right) ,
    \text{ on $\rangset_{\theta_2}$} .
\end{equation}
With this choice we find, using the equalities
(\ref{eq:square_root_rho_equalities})
\begin{equation}
  Q^{-1} \left( \begin{matrix} 0 & \rho \\ -\Arho & 0 \end{matrix} \right) Q
  = \left( \begin{matrix}
        \ii \widetilde{B} & 0 \\
            0 & -\ii \widetilde{B} \end{matrix} \right)
  + \order(-1) ,
  \text{ on $\rangset_{\theta_2}$} ,
\end{equation}
where $\widetilde{B}$ is the square root operator defined in
(\ref{eq:symbol_square_root}).

Next we take the third term in (\ref{eq:transformed_operator2}). It
can be seen easily that
\begin{equation}
  \pdpd{\Arho^{s}}{z} \Arho^{-s}
  = s \pdpd{\Arho}{z} \Arho^{-1} + \order(-1) ,
    \text{ on $\rangset_{\theta_2}$} .
\end{equation}
It follows that
\begin{equation}
  \pdpd{Q^{-1}}{z} Q
    = \left( \begin{matrix}
        0 & \frac{1}{4} \big( \pdpd{\Arho}{z} \Arho^{-1}
                - \pdpd{\rho}{z} \rho^{-1} \big) \\
        \frac{1}{4} \big( \pdpd{\Arho}{z} \Arho^{-1}
                - \pdpd{\rho}{z} \rho^{-1} \big) & 0
      \end{matrix} \right)
  + \order(-1) ,
    \text{ on $\rangset_{\theta_2}$ .}
\end{equation}
Thus, with $Q= Q^{(0)}$, the expression (\ref{eq:transformed_operator2}) is
equal to
\begin{equation} \label{eq:almost_diagonal_operator0}
  \pdpd{}{z}
  - \ii \left( \begin{matrix}
      \widetilde{B} & - \frac{1}{4}\ii \big( \pdpd{\Arho}{z} \Arho^{-1}
                - \pdpd{\rho}{z} \rho^{-1} \big) \\
        - \frac{1}{4}\ii \big( \pdpd{\Arho}{z} \Arho^{-1}
                - \pdpd{\rho}{z} \rho^{-1} \big)
        & - \widetilde{B}
    \end{matrix} \right) + R,
  \text{ on $\rangset_{\theta_2}$}.
\end{equation}
Here $R$ is a pseudodifferential operator of order $1$, that is
of order $-1$ on the set $\rangset_{\theta_2}$.

In expression (\ref{eq:almost_diagonal_operator0}) the highest, first order
part is diagonal, while there are lower order off-diagonal terms. Following
\cite{taylor75} we will remove the off-diagonal terms order by order. To
remove zeroeth order off-diagonal terms we modify $Q$, and set it equal to $Q
= Q^{(1)} = Q^{(0)} ( 1 + K^{(1)} )$. Here $K^{(1)}$ is an operator that
remains to be chosen, is of order $-1$ and of the form $K^{(1)} =
\left(\begin{matrix} 0 & K_{1,2}^{(1)}
\\ K_{2,1}^{(1)} 0 \end{matrix} \right)$. This results in an additional
contribution to (\ref{eq:almost_diagonal_operator0}) of order zero, which is
given to highest, zeroeth order by
\begin{equation} \label{eq:commutator_A_K}
  \left( \begin{matrix}  \ii \widetilde{B} & 0 \\
                           0 & -\ii \widetilde{B} \end{matrix} \right)
                           K^{(1)}
  - K^{(1)}
  \left( \begin{matrix}  \ii \widetilde{B} & 0 \\
                           0 & -\ii \widetilde{B} \end{matrix} \right)
  = \left( \begin{matrix}
  0 & \ii ( \widetilde{B} K^{(1)}_{1,2} + K^{(1)}_{1,2} \widetilde{B}) \\
  -\ii ( \widetilde{B} K^{(1)}_{2,1} + K^{(1)}_{2,1} \widetilde{B}) & 0
       \end{matrix} \right) .
\end{equation}
It follows that there is a symbol $K^{(1)}$ such that this
contribution cancels to zeroeth order off-diagonal contribution on
$\rangset_{\theta_2}$, not changing the first and zeroeth order
diagonal part. By considering further modifications of the form
$Q^{(j+1)} = Q^{(j)}(1+K^{(j+1)})$, $K^{(j+1)}$ of order $-j$, also
the lower order off-diagonal terms can be removed. This proves the
existence of $B_\pm$ satisfying (\ref{eq:symbol_B}).

Next we prove the self adjointness. The operator $Q^{(0)}$ satisfies
\begin{equation}
  \left(\begin{matrix} 1&0 \\ 0&-1 \end{matrix}\right)
    Q^{(0)}{}^* \left(\begin{matrix} 0&-\ii S \\ \ii S&0 \end{matrix}\right)
  =Q^{(0)}{}^{-1}
\end{equation}
It follows that with $Q = Q^{(0)}$ we have
\begin{equation}
  Q^{-1} \left( \begin{matrix} 0 & \rho \\ -\Arho & 0 \end{matrix} \right) Q
  = - \left(\begin{matrix} 1&0 \\ 0&-1 \end{matrix}\right)
    \left[ Q^{-1} \left(
            \begin{matrix} 0 & \rho \\ -\Arho & 0 \end{matrix} \right) Q
            \right]^*
    \left(\begin{matrix} 1&0 \\ 0&-1 \end{matrix}\right)
\end{equation}
Hence with $Q = Q^{(0)}$ the second term
(\ref{eq:transformed_operator2}) is the sum of an anti-selfadjoint
diagonal part and a selfadjoint off-diagonal part (this also follows
from explicit computation). The same property is true for
$\pdpd{Q^{-1}}{z} Q$.  To prove that $B_\pm$ can be chosen selfadjoint
it is sufficient to show that this property is still true when the $Q$
is modified order by order as described above. So suppose the property
holds for $Q^{(j)}$. Then there is a off-diagonal, self-adjoint
$K^{(j+1)}$, with correct highest order term as above. We modify the
definition of $Q^{(j+1)}$ with lower order terms according to
$Q^{(j+1)} = Q^{(j)} \exp(K^{(j+1)})$ (meaning the power series for
$\exp$). Then the microlocal inverse is given by $Q^{(j+1)}{}^{-1} =
\exp(-K^{(j+1)}) Q^{(j)}{}^{-1}$. It is easy to see from the power
series for $\exp$ that
\begin{equation}
  \left(\begin{matrix} 1&0 \\ 0&-1 \end{matrix}\right)
    \exp(-K^{(j+1)})
    \left(\begin{matrix} 1&0 \\ 0&-1 \end{matrix}\right)
  = \exp(K^{(j+1)})
\end{equation}
It follows that the new matrix valued differential operator stays the
sum of an anti-selfadjoint diagonal part and a selfadjoint
off-diagonal part. This shows the selfadjointness.

To compute the second choice of $Q$ and $B_\pm$ we replace $Q$ with $Q D$ where
\begin{equation}
  D = \operatorname{diag}( Q_{11}^{-1} , Q_{12}^{-1} )
\end{equation}
If we denote the second choice of $Q$ by $\widehat{Q}$, it follows easily that
\begin{equation} \label{eq:transformed_operator_hat}
  \widehat{Q}^{-1} \bigg[ \pdpd{}{z}
        - \left( \begin{matrix} 0 & \rho \\ -\Arho & 0 \end{matrix} \right)
    \bigg] \widehat{Q}
  =
  \pdpd{}{z}
    - \ii \left( \begin{matrix}  B_+ & 0 \\
                             0 & B_- \end{matrix} \right)
    - \pdpd{D^{-1}}{z} D + R .
\end{equation}
The operator $D$ is diagonal and the principal symbol of $D$ is given by
$\operatorname{diag}(\rho^{-\frac{1}{2}}a^{\frac{1}{4}},
\linebreak[2]
\rho^{-\frac{1}{2}}a^{\frac{1}{4}})$. This shows the last part of the
lemma.
\qed
\end{pf*}

\section{Proof of Theorem~\ref{th:square_root_correct_sols}%
\label{sec:prove_extrapolation}}

We prove only the $+$ sign, the $-$ sign proceeds in the same
way. Define $v_+, v_-$ by
\begin{equation}
  \left(\begin{matrix} v_+\\v_-\end{matrix}\right)
  = Q^{-1} \left(\begin{matrix} U\\V\end{matrix}\right) .
\end{equation}
Because of (\ref{eq:WF_U_directive}) it follows that $v_-(z_0) \equiv
0$ on $\rangset_{\theta_2} \cap \{ z=z_0 \}$, hence
\begin{equation}
  v_+(z_0) \equiv Q_+^{-1} U(z_0) .
\end{equation}
Since on $J_+(z_0,\theta_2)$ we have $-\tau^{-1} \zeta > 0$, it follows that
$U \equiv Q_+ v_+$ on $J_+(z_0,\theta_2)$.
Equation (\ref{eq:PU_equiv_0}) implies that
\begin{equation}
  \big( \pdpd{}{z} - \ii B_+ \big) v_+ = \widetilde{f}
\end{equation}
for some $\widetilde{f}$ satisfying $\widetilde{f} \equiv 0$ on
$\angleset_{\theta_2}$.  Let $w_+$ be the solution to
\begin{align}
  P_{0,+} w = {}& 0 , &
    w_+|_{z_0} = {}& v_+|_{z_0} .
\end{align}
By the initial condition and a propagation of singularities result, and the fact
that $P_+(v_+ - w_+) \equiv 0$ on $\angleset_{\theta_2}$, it follows that
\begin{equation} \label{eq:v_equiv_w}
  v_+ \equiv w_+ + g ,
\end{equation}
where $g \equiv 0$ on $J_+(z_0,\theta_2)$. Since the symbol of $K_+$ is
in $S^{-\infty}$ outside $J_+(z_0,\theta_2)$, we have $K_+ g \equiv 0$.
Hence
\begin{equation}
  u_+ \equiv K_+ w_+
  \equiv K_+ (v_+ - g)
  \equiv K_+ ( Q_+^{-1} U ) .
\end{equation}
Therefore
\begin{equation}
  Q_+ u_+ \equiv ( K + [Q_+, K] Q_+^{-1} )  U .
\end{equation}
For the commutator term $[Q_+, K] Q_+^{-1}$ it follows from
(\ref{eq:K_suppress_property_1}) and (\ref{eq:K_suppress_property_2}) that
\begin{align}
  [Q_+, K] Q_+^{-1} \sim {}& 0 \text{ on $_+J(z_0,\theta_1)$} ,
    \label{eq:comm_suppress_property_1} \\
  [Q_+, K] Q_+^{-1} \in {}& \Op S^{-\infty}
                \text{ outside $J_+(z_0,\theta_2)$,
        if $z-z_0 > \delta$} .
    \label{eq:comm_suppress_property_2}
\end{align}
This completes the proof of the theorem.


\def\cprime{$'$}

\end{document}